\documentclass[11pt]{article}

\usepackage{amsmath,amssymb}
\usepackage{latexsym}
\usepackage{graphicx}
\usepackage{bm}
\newcounter{stmtcount}
\newcounter{algocount}
{\catcode`\;=\active%
\gdef\pointvirgule{\catcode`\;=\active%
\def;{\unskip\kern 2pt\string;~~~\ignorespaces}}}
\newenvironment{ualgorithm}[1]{\pointvirgule\setcounter{stmtcount}{0}
\def\N>{\stepcounter{stmtcount}\thestmtcount.\>}\def\B>{\hskip0.1em\vrule\>}
\def\NN>{\stepcounter{stmtcount}\ifnum\thestmtcount<10\phantom{0}\fi
\thestatement.\>}\def\[##1]{{\bf##1}}
\def\[##1]{{\bf##1}}
\refstepcounter{algocount}
\begin{tabbing}
++\=++\=++\=++\=++\=++\=++\=++\=++\=++\=++\=++\=++\=++\kill
\indent \underline{{\sc Algorithm} \thealgocount}: #1\\}%
{\end{tabbing}}
\begin{document}

\markboth{N. S. Hoang and A. G. Ramm}{Solving ill-conditioned las by DSM}

\title{
Solving ill-conditioned linear algebraic systems by\break
the dynamical systems method (DSM)}

\author{N. S. Hoang\footnotemark[1] ~and
	  A. G. Ramm\footnotemark[2]\\
Mathematics Department, Kansas State University,\\
Manhattan, KS 66506, USA\\
email: \{nguyenhs, ramm\}@math.ksu.edu}%~~~}

\maketitle

\begin{abstract}
An iterative scheme for the Dynamical Systems Method (DSM) is given such that 
one does not have to solve the Cauchy problem occuring in the application 
of the DSM for solving ill-conditioned linear algebraic systems. 
The novelty of the algorithm is that the algorithm does not have to find 
the regularization parameter $a$ by solving a nonlinear equation. 
Numerical experiments show that DSM competes favorably with the Variational Regularization.

%\begin{keywords}
%{\small ill-posed problems, ill-conditioned linear algebraic systems, dynamical systems method (DSM).}
%\end{keywords}

{\bf Keywords}:
{\small ill-posed problems, ill-conditioned linear algebraic systems, dynamical systems method (DSM).}

{\bf AMS subject classification}:
{\small 65F10, 65F22.}
%\begin{classcode}
%{\small 65F10, 65F22.}
%\end{classcode}
\end{abstract}

\section{Introduction}

Consider a linear operator equation of the form
\begin{equation}
\label{linear}
F(u)=Au-f=0,\quad u\in H,
\end{equation}
where $H$ is a Hilbert space and $A$ is a linear operator in $H$ which is 
not necessarily bounded but closed and densely defined.
To solve this equation we apply the Dynamical Systems Method (DSM) introduced in \cite{Ramm}:
\begin{equation}
\label{cauchy}
u'=-u+(T+a(s))^{-1}A^{*}f,\quad u(0)=u_0,
\end{equation}
where $T:=A^{*}A$ and $a(t)>0$ is a nonincreasing function such that $a(t)\to 0$ as $t\to\infty$. 
The 
unique solution to \eqref{cauchy} is given by
\begin{equation}
\label{integral}
u(t)=u_0e^{-t}+e^{-t}\int_0^te^{s}(T+a(s))^{-1}A^{*}fds.
\end{equation}

The DSM consists of solving problem \eqref{cauchy} with a chosen $a(t)$ 
and $u_0$ and finding a stopping time $t_\delta$ so that  $u(t_{\delta})$ 
approximates the solution $y$ to problem \eqref{linear} of minimal norm. 
Different choices of $a(t)$ generate different methods of solving equation \eqref{linear}. 
These methods have different accuracy and different computation time. 
Thus, in order to get an efficient implementation of the DSM, we need to study the choice of $a(t)$ 
and of the stopping time $t_\delta$. Since the solution to \eqref{linear} can be presented in 
the form of an integral, the question arises: how can one compute the integral efficiently? 
The integrand of the solution is used also in the Variational Regularization (VR) method. 
The choice of the stopping time $t_\delta$ will be done by a discrepancy-type principle for DSM. 
However, choosing $a(t)$ so that the method will be accurate and the computation time is small
 is not a trivial task.

This paper deals with the following questions:
\begin{enumerate}
\item How can one choose $a(t)$ so that the DSM is fast and accurate?
\item Does the DSM compete favorably with the VR in computation time?
\item Is the DSM comparable with the VR in accuracy?
\end{enumerate}

\section{Construction of method}

\subsection{An iterative scheme}
\label{iterativesec}

Let us discuss a choice of $a(t)$ which allows one to solve problem \eqref{cauchy}
 or to calculate the integral \eqref{integral} without using any numerical method for solving 
initial-value problem for ordinary differential equations (ODE). 
In fact, using a monotonically decreasing $a(t)$ with one of the best numerical methods for nonstiff ODE,
 such as DOPRI45, is more expensive computationally than using a step function $\tilde{a}(t)$,
 approximating $a(t)$, but 
brings no improvement in the accuracy of the solution to our problems compared to 
the numerical solution of our problems given in Section \ref{doprivsstep}.  

Necessary conditions for the function $a(t)$ are: $a(s)$ is a nonincreasing function 
and $\lim_{s\to\infty} a(s)=0$ (see \cite{Ramm}). Thus, our choice of $a(t)$ must satisfy these conditions. 
Consider a step function $\tilde{a}(t)$, approximating $a(t)$, defined as follows:
$$
\tilde{a}(t)=a_n,\quad t_n\leq t<t_{n+1},
$$
the number $t_n$ are chosen later.
For this $\tilde{a}(t)$, $u_n=u(t_n)$ can be computed by the formula:
\begin{align*}
u_n=u_0e^{-t_n}+e^{-t_n}\sum_{i=1}^n(e^{t_i}-e^{t_{i-1}})(T+a_{i-1})^{-1}A^*f_\delta.
\end{align*}
This leads to the following iterative formula:
\begin{equation}
\label{iterfor}
u_{n+1}=e^{-h_n}u_{n}+(1-e^{-h_n})\big{(}T+a_n\big{)}^{-1}A^*f_\delta,\quad h_n=t_{n+1}-t_n.
\end{equation}
Thus, $u_{n}$ can be obtained iteratively if $u_0$\,, $a(t)$ and $t_{n}$ are known. 

The questions are:
\begin{enumerate}
\item{For a given $a(t)$, how can we choose $t_n$ or $h_n$ so that the DSM works efficiently?}
\item{With $a_n=a(t_n)$ where $a(t)$ is a continuous function, does the iterative scheme 
compete favorably with the DSM version in which $u(t)$ is solved by some numerical methods 
such as Runge-Kutta methods using $a(t)$?}
\end{enumerate}

In our experiments, $a_n=a(t_n)$ where $a(t)=\frac{a_0}{1+t}$ where $a_0>0$ is a constant
which will be chosen later, 
as suggested in \cite{Ramm}, with $t_n$ chosen so that $t_{n+1}-t_{n}=h_n$, $h_n=q^n$, 
where $1\leq q\leq2$. For this choice, if $q>1$ then the solution  $u_n$ at the $n$-th 
step depends mainly on $\big{(}T+ a_n\big{)}^{-1}A^*f_\delta$  since $e^{-h_n}$ 
is very small when $n$ is large. 

Note that $a_n$ decays exponentially fast when $n\to\infty$ if $q>1$. 
A question arises: how does one choose $q$ so that the method is fast and accurate? 
This question will be discussed in Section 3.

ALGORITHM~\ref{DSM}
demonstrates the use of the iterative formula \eqref{iterfor} and a relaxed discrepancy
principle described below for finding $u$ given $a_0$, $A$, $f_\delta$ and $\delta$.

In order to improve the speed of the algorithm, we use a relaxed discrepancy principle: 
at each iteration one checks if 
\begin{equation}
\label{discrep1}
0.9\delta \leq \|Au_n-f_\delta\| \leq 1.001\delta.
\end{equation}
As we shall see later, $a_0$ is chosen so that the condition \eqref{a0} (see below) is satisfied. 
Thus, if $u_0=T_{a_0}^{-1}A^*f_\delta$, where $T_{a}:=T+a$, then $\delta < \|Au_0-f_\delta\|$.
Let $t_n$ be the first time such that 
$\|Au_n-f_\delta\|\leq 1.001\delta $. If \eqref{discrep1} is satisfied, then one 
stops calculations.
If $\|Au_n-f_\delta\|<0.9\delta $, then one takes a smaller step-size and 
recomputes $u_n$.
If this happens, we do not increase $h_n$, that is, we do not multiply 
$h_n$ by $q$ in the following steps. 
One repeats this procedure until condition \eqref{discrep1} is satisfied. 

\begin{center}
\begin{tabular}[b]{c}\vbox{%hide this line to get no box
\begin{ualgorithm}{DSM$(A,f_\delta,\delta)$}
$q:= 2$;\\
$g_\delta:=A^*f_\delta$;  $T:=A^*A$;       \\
$itermax := 30$; $u=(T+a_0)^{-1}g_\delta$;\\
$i:=0$; $t=1$; $h:=1$; $halve:=0$;\\
\[while] $(1.001\delta <\|Au-f_\delta\|)$ \[and] $(i<itermax)$ \[do]  \\
\B> $i:=i+1$; $t=t+h$; $a=a_0/t$;\\
\B> $v := (T+aI)^{-1} g_\delta$; \\
\B> $\tilde{u}=e^{-h}u+(1-e^{-h})v$; \\
\B> \[if] $0.9\delta < \|A\tilde{u}-f_\delta\|$ \[then] \\
\B> \> $u := \tilde{u}$; \\
\B> \> \[if] $halve=0$ \[then] $h := hq$; \[end]; \\
\B> \[elseif]\\
\B> \> $t:=t-h$; $h:= h/2$; $halve = 1$;\\
\B> \[endif]\\
\[endwhile]
\label{DSM}
\end{ualgorithm}
} \\ \end{tabular}
\end{center}

In order to improve the speed of the algorithm, we use a relaxed discrepancy principle: 
at each iteration one checks if 
\begin{equation}
\label{discrep1}
0.9\delta \leq \|Au_n-f_\delta\| \leq 1.001\delta.
\end{equation}
As we shall see later, $a_0$ is chosen so that the condition \eqref{a0} (see below) is satisfied. 
Thus, if $u_0=T_{a_0}^{-1}A^*f_\delta$, where $T_{a}:=T+a$, then $\delta < \|Au_0-f_\delta\|$.
Let $t_n$ be the first time such that 
$\|Au_n-f_\delta\|\leq 1.001\delta $. If \eqref{discrep1} is satisfied, then one 
stops calculations.
If $\|Au_n-f_\delta\|<0.9\delta $, then one takes a smaller step-size and 
recomputes $u_n$.
If this happens, we do not increase $h_n$, that is, we do not multiply 
$h_n$ by $q$ in the following steps. 
One repeats this procedure until condition \eqref{discrep1} is satisfied. 

\subsection{On the choice of $a_0$}

From numerical experiments with ill-conditioned linear algebraic systems (las) of the form
$Au=f_\delta$, it follows that
 the regularization parameter $a_M$, obtained from the discrepancy
 principle $\|Au_{a_M}-f_\delta\|=\delta$, where $u_{a_M}=T_{a_M}^{-1}A^*f_\delta$, is often 
close to the optimal value $a_{op}$, i.e., the value minimizing the quantity:
\begin{equation*}
\|u_{a_{op}}-y\| = \inf_{a}\|u_a -y\|,\quad
u_a = T_a^{-1}Af_\delta.
\end{equation*}
The letter $M$ in $a_M$ stands for
Morozov, who suggested to choose $c=1$ in the disrepancy principle. 

If $a_0$ is chosen smaller than $a_{op}$, the method may converge poorly.
Since $a_M$ is close to $a_{op}$, only those $a$ for which $\|AT_a^{-1}A^*f_\delta-f_\delta\|=c\delta$
 with $c$ 'close' to $1$ yield accurate approximations to the solution $y$. 
Also, if $a_0$ is chosen much greater than $a_{op}$, then the information obtained from the starting 
steps of the iterative process \eqref{iterfor}
is not valuable because when $a_0$ is far from $a_{op}$, the error $\|u_{a_0}-y\|$ is much bigger than
$\|u_{a_{op}}-y\|$.
If $a_0$ is much bigger than $a_{op}$, a lot of time will be spent until $a(t_n)$ becomes close to $a_{op}$. 
In order to increase the speed of computation, $a_0$ should be chosen so that
 it is close to $a_{op}$ and greater than $a_{op}$. 
Since $a_{op}$ is not known and is often close to $a_M$, we choose $a_0$ from the condition: 
\begin{equation}
\label{a0}
\delta <\|Au_{a_0}-f_\delta\|<2\delta.
\end{equation}
For this choice, $a_0$ is 'close' to and greater than $a_M$. 
Since there are many $a_0$ satisfying this condition, it is not difficult to find one of them.

In the implementation of the VR using discrepancy principle with Morozov's suggestion $c=1$, 
if one wants to use the Newton method for finding the regularization parameter, 
one also has to choose the starting value $a_0$ so that the iteration process converges, 
because the Newton method, in general, converges only locally. 
If this value is close to and greater than $a_{M}$, it can also be used as the initial value of 
$a_0=a(t)|_{t=0}$ for the DSM. 

In our numerical experiments, with a guess $a_0=\frac{1}{3}\max\lambda_i(A^*A)\delta_{rel}$ 
for $a(0)$, we find $a_0$ such that 
$\delta <\|Au_{a_0}-f_\delta\|<2\delta$. Here, $\delta_{rel}$ stands for the relative error,
 i.e., $\delta_{rel}=\frac{\delta}{\|f\|}$.
The factor $\frac{1}{3}$ is introduced here in order to reduce the cost for finding $a_0$,
because $a_0$, which satisfies \eqref{a0}, is often less than $\max \lambda_i(A^*A)\delta_{rel}$.
The idea for this choice is based on the fact that the spectrum of the matrix $\frac{1}{\max \lambda_i(A^*A)}A^*A$ 
is contained in $[0,1]$. 

Note that ones has
$$
a_M\leq \frac{\delta \|A\|^2}{\|f_\delta\|-\delta}.
$$
Indeed,
\begin{align*}
\|f_\delta\|-\delta &= \|f_\delta\|-\|Au_{a_M}-f_\delta\|\leq \|Au_{a_M}\|.
\end{align*}
Since $A^*Au_{a_M}+a_{M}u_{a_M}=A^*f_\delta$, one has $a_MAu_{a_M}= AA^*(f_\delta- Au_{a_M})$. 
Thus,
$$
\|f_\delta\|-\delta\|\leq \|Au_{a_M}\|=\frac{1}{a_M}\|AA^*(f_\delta- Au_{a_M})\|\leq \frac{\|A\|^2}{a_M}\delta.
$$
Similar estimate one can find in \cite[p. 53]{Kirsch}, where $a_0=\frac{\delta \|A\|^2}{\|f_\delta\|-\delta}$ 
is suggested as a starting value for Newton's method to determine $a_M$ on the basis that
it is an upper bound for $a_M$. 
Note that $\frac{\delta \|A\|^2}{\|f_\delta\|-\delta}\approx\delta_{rel}\|A\|^2=\max \lambda_i(A^*A)\delta_{rel}$.
However, in practice 
Newton's method does not necessarily converge with this starting value. If this happens,
a smaller starting value $a_1:=\frac{a_0}{2}$ is used to restart the Newton's method.

In general, our initial choice for $a_0$ may not satisfy \eqref{a0}. 
Iterations for finding $a_0$ to satisfy \eqref{a0} are done as follows:

\begin{enumerate}
\item{
If  $\frac{\|Au_{a_0} -f_\delta\|}{\delta}=c>3$, then one takes $a_1:=\frac{a_0}{2(c-1)}$ as the next guess
 and checks if the condition \eqref{a0} is satisfied. If $2<c\leq 3$ then one takes $a_1:=a_0/3$.}
\item{
If  $\frac{\|Au_{a_0} -f_\delta\|}{\delta}=c<1$, then $a_1:=3a_0$ is used as the next guess.}
\item{
After $a_0$ is updated, one checks if \eqref{a0} is satisfied. 
If \eqref{a0} is not satisfied, one repeats steps 1 and 2 until 
one finds $a_0$ satisfying condition \eqref{a0} (see ALGORITHM~\ref{itera0}).} 
\end{enumerate}

\begin{center}
\begin{tabular}[b]{c}\vbox{%hide this line to get no box
\begin{ualgorithm}{find-$a_0$}
$a_0:=\frac{1}{3}\|A\|^2\delta_{rel}$;\\
$c:=\|A u_{a_0}-f_\delta\|/\delta$;\\
\[while] $(2 < c)$ \[or] $(c<1)$ \[do]                        \\
\B> \[if] $3<c$ \[then] \\
\B> \> $a_0:=0.5a_0/(c-1)$;\\
\B> \[elseif] $(2 < c\leq 3)$ \[then]\\
\B> \> $a_0 := a_0/3$; \\
\B> \[else]\\
\B> \> $a_0:=3a_0;$\\
\B> \[end]\\
\B>  $u_{a_0}:=(A^*A+a_0)^{-1}A^*f_\delta$;\\
\B> $c:=\|A u_{a_0}-f_\delta\|/\delta$;\\
\[endwhile]
\label{itera0}
\end{ualgorithm}
} \\ \end{tabular}
\end{center}

The above strategy is based on the fact that the function 
$$
\phi(a)=\|A(T+a)^{-1}A^*f_\delta-f_\delta\|
$$
is a monotonically decreasing function of $a$, $a>0$. 
In looking for $a_0$, satisfying \eqref{a0}, when our guess
$a_0\gg a_M>0$ or $\|Au_{a_0}-f_\delta\|\gg \delta$, one uses an approximation
\begin{align*}
\phi(x)&\approx \phi(a_0) + (x-a_0)\frac{\phi(a_0)-\phi(a_M)}{a_0-a_M}\\
&\approx \phi(a_0) + (x-a_0)\frac{\phi(a_0)-\phi(a_M)}{a_0}=:\varphi(x).
\end{align*}
Note that $\phi(a_0)$ and $a_0$ are known. We are looking for $x$ such that 
$\delta<\varphi(x)<2\delta$.
Thus, if $a_1$ is such that $\delta<\varphi(a_1)<2\delta$ and if $2\delta<\phi(a_0)$, then
$$
(\phi(a_0)-2\delta)\frac{a_0}{\phi(a_0)-\delta}<a_0-a_1<(\phi(a_0)-\delta)\frac{a_0}{\phi(a_0)-\delta}.
$$
Hence, we choose $a_1$ such that
$$
a_0-a_1=(\phi(a_0)-1.5\delta)\frac{a_0}{\phi(a_0)-\delta},
$$
so
$$
a_1=a_0\frac{0.5\delta}{\phi(a_0)-\delta}.
$$
Although this is a very rough approximation, it works well in practice. 
It often takes 1 to 3 steps to get an $a_0$ satisfying \eqref{a0}. 
That is why we have a factor $\frac{0.5}{c-1}$ in the first case. 
Overall, it is easier to look for $a_0$ satisfying \eqref{a0} than 
to look for $a_0$ for which the Newton's method converges.
Indeed, the Newton's scheme for solving $a_M$ does not necessarily
converge with $a_0$ found from condition \eqref{a0}.

\section{Numerical experiments}

In this section, we compare DSM with VR$_{i}$ and VR$_{n}$.
In all methods, we begin with the guess $a_0=\frac{1}{3}\|A\|^2\delta_{rel}$
and use the ALGORITHM~\ref{itera0} to find $a_0$ satisfying condition~\eqref{a0}.
In our experiments, the computation cost for this step is very low.
Indeed, it only takes 1 or 2 iterations to get $a_0$.
By VR$_i$ we denote the VR obtained by using $a=a_0$, 
the intial value for $a(t)$ in DSM, and by VR$_{n}$ we denote the VR with $a=a_M$, found from 
the VR discrepancy principle with $c=1$ by using Quasi-Newton's method with the initial guess $a=a_0$.
Quasi-Newton's method is chosen instead of Newton's method in order to reduce the computation cost. 
In all experiments we compare these methods in accuracy and with respect to the parameter 
$N_{linsol}$, which is the number of times 
for solving the linear system $T_{a}u=A^*f_\delta$ for $u$. 
Note that solving these linear systems is the main cost in these methods.

In this section, besides comparing the DSM with the VR for linear algebraic systems with Hilbert matrices, 
we also carry out experiments with other linear algebraic systems given in 
the {\it Regularization package} in \cite{Hansen1}. 
These linear systems are obtained as a part of numerical solutions to some integral equations. 
Here, we only focus on the numerical methods for solving linear algebraic systems, 
not on solving these integral equations. 
Therefore, we use these linear algebraic systems to test our methods
for solving stably these systems.

\subsection{Linear algebraic systems with Hilbert matrices}
Consider a linear algebraic system
\begin{align}
\label{hihi}
H_{n}u=f_{\delta},
\end{align}
where

$$
f_{\delta}
=f
+e,
\quad
f=H_{n}x, 
\quad
H_{n}= 
\begin{bmatrix}
1&\frac{1}{2}&\cdots&\frac{1}{n}\\
\frac{1}{2}&\frac{1}{3}&\cdots&\frac{1}{n+1}\\
\vdots&\vdots&\ddots&\vdots\\
\frac{1}{n}&\frac{1}{n+1}&\cdots&\frac{1}{2n-1}
\end{bmatrix},
$$
and $e\in \mathbb{R}^n$ is a random normally distributed vector such that $\|e\|_2 \leq \delta_{rel}\|f\|_2$. 
The Hilbert matrix $H_n$ is well-known for having a very large condition number when $n$ is large. 
%In fact, the condition number
%of $H_n$ is estimated by the formula $\cond(H_n)=O(e^{3.5255 n}/\sqrt{n})$.
If $n$ is sufficiently large, the system is severely ill-conditioned.

\subsubsection{The condition numbers of Hilbert matrices}

It is impossible to calculate the condition number of $H_n$ by
computing the ratio of the largest and the smallest eigenvalues of $H_n$ 
because for large $n$ the smallest eigenvalue of $H_n$ is smaller than 
$10^{-16}$. 
Note that singular values of $H_n$ are its eigenvalues since $H_n$ is selfadjoint and positive definite.
Due to the limitation of machine precision, every value smaller than 
$10^{-16}$ is understood as 0. 
That is why if we use the function {\it cond} provided by MATLAB, the condition number of $H_n$ 
for $n\geq20$ is about $10^{16}\times \max |\lambda_i(H_n)|$.
Since the largest eigenvalue of $H_n$ grows very slowly, the condition numbers 
of $H_n$
for $n\geq 20$ are all about $10^{20}$,
while, in fact, the condition number of $H_{100}$ 
computed by the formula, given below, is about $10^{150}$ (see 
Table~\ref{Hilbcon}). 
In general, computing condition numbers of strongly ill-conditioned matrices is an open problem. 
The function {\it cond}, provided by MATLAB, is not always reliable
for computing the condition number of ill-condition matrices. 
Fortunately, there is an analytic formula for the inverse of $H_n$. 
Indeed, one has (see \cite{choi})
$H_n^{-1}=(h_{ij})_{i,j=1}^n$, where 
$$
h_{ij}=(-1)^{i+j}(i+j-1)\dbinom{n+i-1}{n-j}\dbinom{n+j-1}{n-i}\dbinom{i+j-2}{i-1}^2.
$$
Thus, the condition number of the Hilbert matrix can be computed 
by the formula:
$$
cond(H_n)=\|H_n\|\|H_n^{-1}\|.
$$
Here $cond(H_n)$ stands for the condition number of the Hilbert matrix $H_n$ and $\|H_n\|$ and $\|H_n^{-1}\|$ 
are the largest eigenvalues of $H_n$ and $H_n^{-1}$, respectively. Although MATLAB
can not compute values less than $10^{-16}$, it can compute values up to $10^{200}$. 
Therefore, it can compute $\|H_n^{-1}\|$ for $n$ up to 120.
In MATLAB, the matrices $H_n$ and $H_n^{-1}$ can be obtained by 
the syntax: $H_n=\text{hilb}(n)$ and $H_n^{-1}=\text{invhilb}(n)$, respectively. 

The condition numbers of Hilbert matrices, computed by the above formula, are given in Table~\ref{Hilbcon}.

\begin{table}[h] 
\caption{The condition number of Hilbert matrices.}
\label{Hilbcon}
\centering
\small
\begin{tabular}{@{  }r@{\hspace{2mm}}@{\hspace{2mm}}|c@{\hspace{2mm}}
c@{\hspace{2mm}}c@{\hspace{2mm}}c@{\hspace{2mm}}c@{\hspace{2mm}}c@{\hspace{2mm}}} 
$n$ &20&40&60&80&100&120\\
\hline
  $cond(H_n)$
  &$2.5\times10^{28}$
  &$7.7\times10^{58}$
  &$2.7\times10^{89}$
  &$9.9\times10^{119}$
  &$3.8\times10^{150}$
  &$1.5\times10^{181}$\\
\end{tabular}
\end{table}
From Table~\ref{Hilbcon} one can see that the computed condition numbers of the Hilbert matrix grow very fast as $n$ grows.
%agree 
%with the formula $\cond(H_n)=O(e^{3.5255 n}/\sqrt{n})$. Note that $e^{3.5255}=10^{1.5311}$.

\subsubsection{Continuous $a(t)$ compared to the step function $\tilde{a}(t)$}
\label{doprivsstep}

In this section, we compare the DSM, which is implemented by solving the Cauchy problem \eqref{cauchy} with $a(t)$, 
and the iterative DSM implemented with $\tilde{a}(t)$ approximating $a(t)$ as
described in Section~\ref{iterativesec}. Both of them use the same $a_0$ which is 
found by ALGORITHM~\ref{itera0}. 
The DSM using a numerical method to solve the Cauchy problem is implemented as follows:

\begin{enumerate}
\item{
One uses the DOPRI45 method which is an embedded pair consisting of a Runge-Kutta (RK) 
method of order 5 and another RK method of order 4 which is used to estimate the error in order to 
control the step sizes. The DOPRI45 is an explicit method which requires 
6 right-hand side function evaluations at each step. 
Details about the coefficients and variable step size strategy can be found in \cite{BURRAGE95, Hairer}. 
Using a variable step size helps to choose the best step sizes and improves the speed.
}
\item{
In solving \eqref{cauchy}, at the end of each step, one always checks the stopping rule, 
based on the discrepancy principle 
$$
0.9\leq\|Au_{\delta}(t)-f_{\delta}\|\leq 1.001\delta.
$$ 
If this condition is satisfied, one stops and takes the solution at the final step $u(t_n)$ 
as the solution to the linear algebraic system.
}
\end{enumerate}

\begin{table}[h] 
\caption{Numerical results for Hilbert matrices for $\delta_{rel}=0.01$, $n=100$.}
\label{convsdis}
\centering
\small
\begin{tabular}{@{  }c@{\hspace{2mm}}|c@{\hspace{2mm}}
@{\hspace{2mm}}c@{\hspace{2mm}}|c@{\hspace{2mm}}c@{\hspace{2mm}}
%|c@{\hspace{2mm}}c@{\hspace{2mm}}
|c@{\hspace{2mm}}c@{\hspace{2mm}}|c@{\hspace{2mm}}r@{\hspace{2mm}}l@{}} 
\hline
&
\multicolumn{2}{c|}{DSM}
%&\multicolumn{2}{c|}{DSM1}
&\multicolumn{2}{c|}{DSM($q=1$)}&\multicolumn{2}{c|}{DSM-DOPRI45}\\
$n$&
$N_{\text{linsol}}$&$\frac{\|u_\delta-y\|_{2}}{\|y\|_2}$&
%$N_{\text{linsol}}$&$\frac{\|u_\delta-y\|_{2}}{\|y\|_2}$&
$N_{\text{linsol}}$&$\frac{\|u_\delta-y\|_{2}}{\|y\|_2}$&
$N_{\text{linsol}}$&$\frac{\|u_\delta-y\|_{2}}{\|y\|_2}$\\
\hline
10    &5    &0.1222   &10    &0.1195  &205    &0.1223\\
20    &5    &0.1373    &7    &0.1537  &145    &0.1584\\
30    &7    &0.0945   &20    &0.1180  &313    &0.1197\\
40    &5    &0.2174    &7    &0.2278  &151    &0.2290\\
50    &6    &0.1620   &14    &0.1609  &247    &0.1609\\
60    &6    &0.1456   &16    &0.1478  &253    &0.1480\\
70    &6    &0.1436   &13    &0.1543  &229    &0.1554\\
80    &6    &0.1778   &10    &0.1969  &181    &0.1963\\
90    &6    &0.1531   &13    &0.1535  &307    &0.1547\\
100   &7    &0.1400   &23    &0.1522  &355    &0.1481\\
\hline 
\end{tabular}
\end{table}
The DSM version implemented with the DOPRI45 method is denoted DSM-DOPRI45 while the other 
iterative version of DSM is denoted just by DSM. 

Table~\ref{convsdis} presents the numerical results with Hilbert matrices
$H_n$ obtained by two versions of the DSM for $n=10,20,...,100$,
$\delta_{rel}=0.01$, $x=(x_1,...,x_n)^T$,
$x_i=\sqrt{2\frac{i-1}{100}\pi}$.  
From Table~\ref{convsdis}, as well as other
numerical experiments, we found out that the accuracy obtained by the
DSM-DOPRI45 is worse than that of the iterative DSM.  Moreover, the
computation time for the DSM-DOPRI45 is much greater than that for the
iterative DSM.  Also, using $h=$const or $q=1$ does not give more accurate
solutions while requires more computation time.

The conclusion from this experiment as well as from other experiments is that the DSM with $q=2$ 
is much faster and often gives better results than the DSM with $q=1$ and the DSM-DOPRI45. 
Therefore, we choose the iterative DSM with $q=2$ to compare with the VR$_{n}$ method.

\subsubsection{DSM compared to VR}

\begin{figure}[!h!tb]
\centerline{%
\includegraphics[scale=0.45]{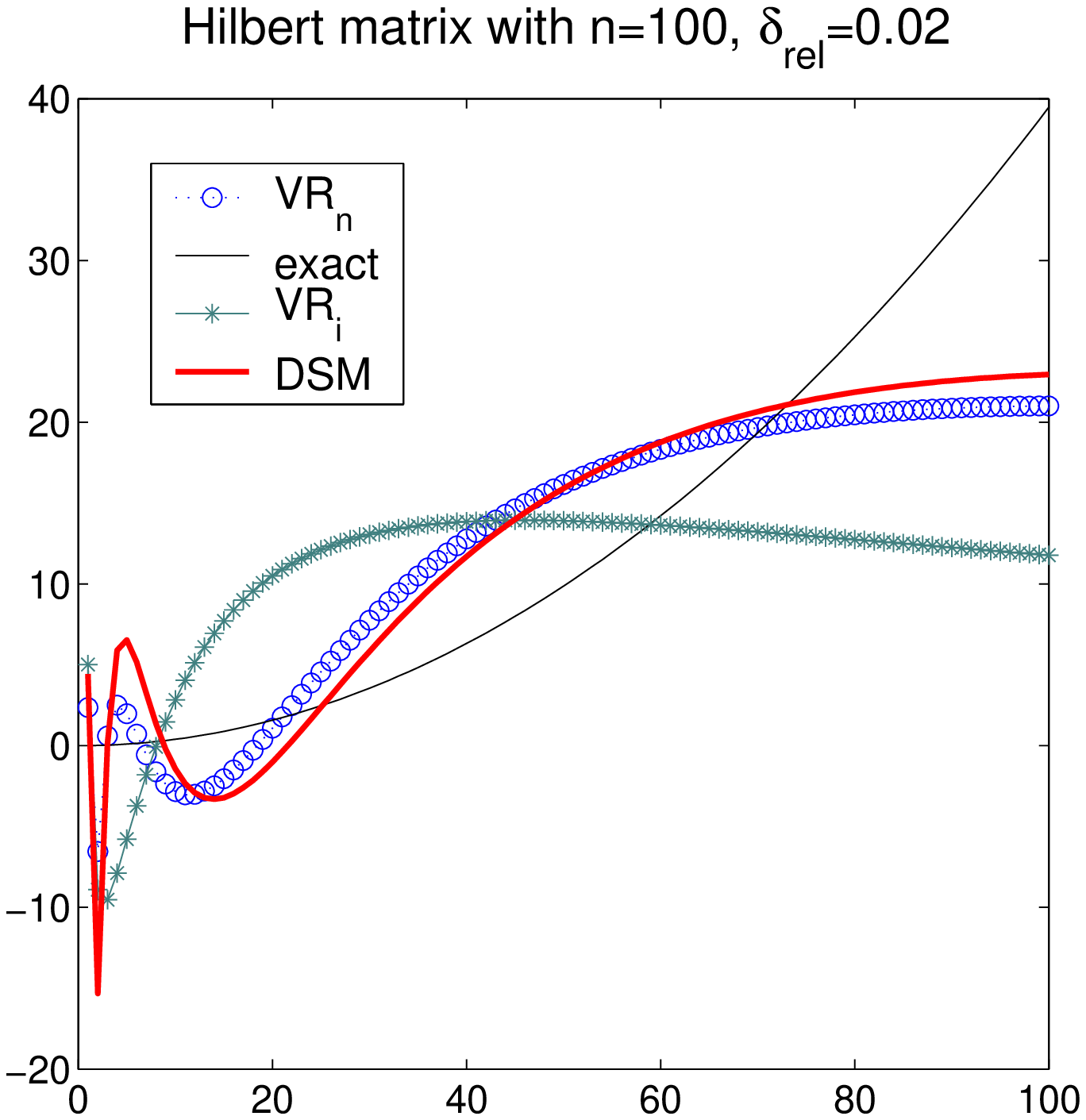}\quad\includegraphics[scale=0.45]{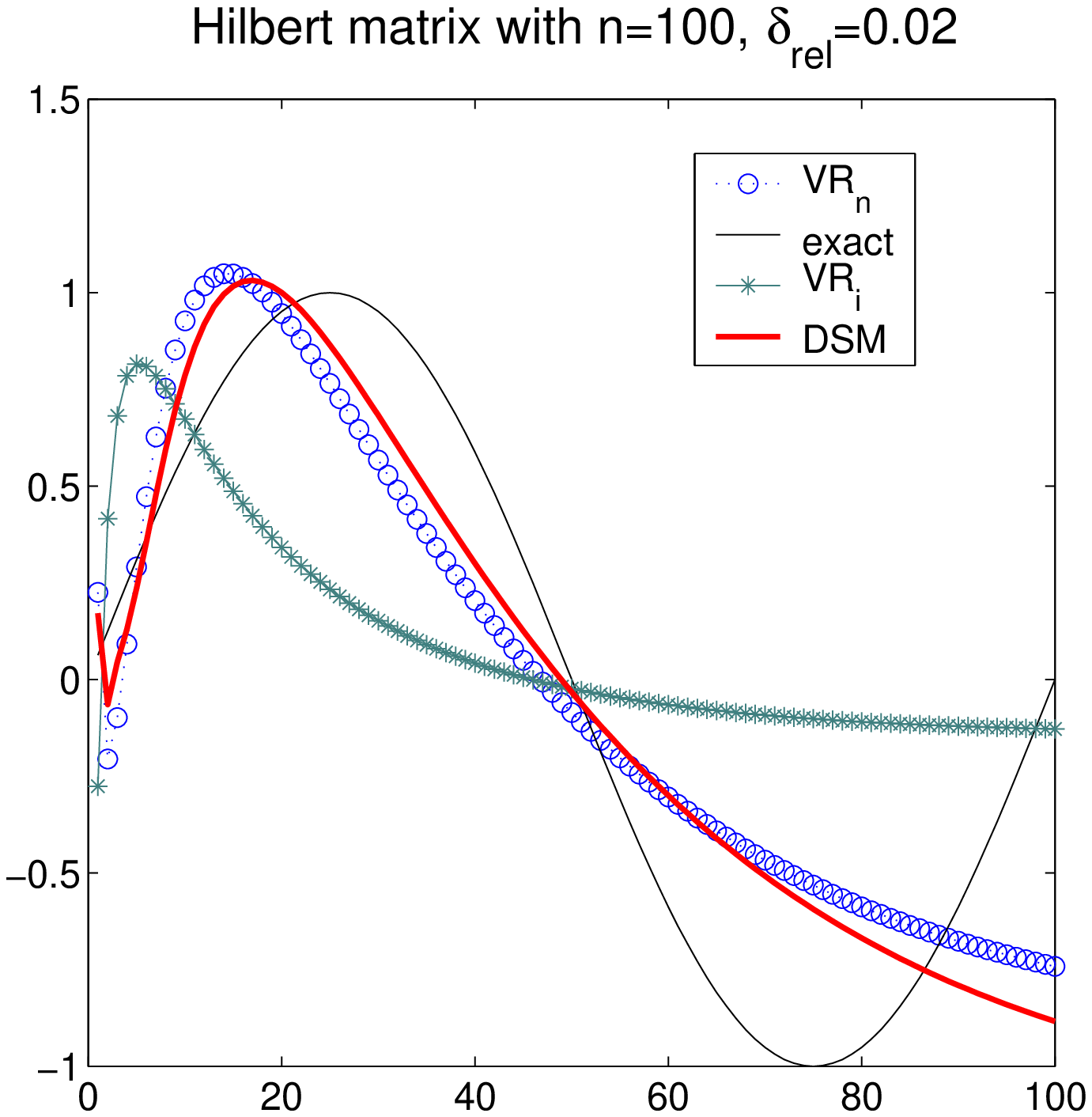}}%{data2.eps}}
\caption{Plots of solutions obtained by the DSM and VR with the exact solution $x$, $x=(x_i)_{i=1}^{100}$ 
when $x_i=(2\frac{i-1}{100}\pi)^2$ (left) 
and $x_i=\sin(2\frac{i-1}{100}\pi)$ (right) with $\delta_{rel} = 0.02$.}
\label{fig1}
\end{figure}

In this section, we test three methods: the DSM, the VR$_{i}$ and the VR$_{n}$ on linear algebraic systems 
with Hilbert matrices. The first linear system is obtained by taking $H_{100}$ and 
$x=(x_1,...,x_{100})^{T}$, where $x_i=(\frac{i-1}{100})^2$. For the second problem we just change $x_i$ to
$x_i=\sin(2\frac{i-1}{100}\pi)$.
Numerical results for these systems are shown in Figure~\ref{fig1}.
%Numerical experiments show that DSM competes favorably with the VR (see figures).

Looking at Figure~\ref{fig1}, one can see that with the same guess $a_0$, both the VR$_{n}$ and
DSM give better results than those of VR$_{i}$. As it can be seen from Figure~\ref{fig1}, 
the numerical solutions obtained by the DSM in these tests are slightly more accurate than those of the VR$_{n}$. 

Table~\ref{hilbert2} presents results with Hilbert matrices $H_n$ for
$n=10,20,...,100$, $\delta_{rel}=0.01$, $x=(x_1,...,x_n)^T$,
$x_i=\sqrt{2\frac{i-1}{100}\pi}$.  Looking at this Table it is clear that
the results obtained by the DSM are slightly more accurate than those by
the VR$_{n}$ even in the cases when the VR$_{n}$ requires much more work
than the DSM.  In this example, we can conclude that the DSM is better
than the VR$_{n}$ in both accuracy and time of computation.

\begin{table}[h] 
\caption{Numerical results for Hilbert matrix $H_n$ for $\delta_{rel}=0.01$, $n=10,20,...,100$.}
%Numerical solution for t = pi*t;x = sqrt(t)';
%
\label{hilbert2}
\centering
\small
\begin{tabular}{@{  }c@{\hspace{2mm}}|
@{\hspace{2mm}}c@{\hspace{2mm}}c@{\hspace{2mm}}|c@{\hspace{2mm}}c@{\hspace{2mm}}
|c@{\hspace{2mm}}c@{\hspace{2mm}}|c@{\hspace{2mm}}r@{\hspace{2mm}}l@{}} 
\hline
&\multicolumn{2}{c|}{DSM}
&\multicolumn{2}{c|}{VR$_{i}$}&\multicolumn{2}{c|}{VR$_{n}$}\\
n&
$N_{\text{linsol}}$&$\frac{\|u_\delta-y\|_{2}}{\|y\|_2}$&
$N_{\text{linsol}}$&$\frac{\|u_\delta-y\|_{2}}{\|y\|_2}$&
$N_{\text{linsol}}$&$\frac{\|u_\delta-y\|_{2}}{\|y\|_2}$\\
\hline
10    &4    &0.2368    &1    &0.3294    &7    &0.2534\\
20    &5    &0.1638    &1    &0.3194    &7    &0.1765\\
30    &5    &0.1694    &1    &0.3372   &11    &0.1699\\
40    &5    &0.1984    &1    &0.3398    &8    &0.2074\\
50    &6    &0.1566    &1    &0.3345    &7    &0.1865\\
60    &5    &0.1890    &1    &0.3425    &8    &0.1980\\
70    &7    &0.1449    &1    &0.3393   &11    &0.1450\\
80    &7    &0.1217    &1    &0.3480    &8    &0.1501\\
90    &7    &0.1259    &1    &0.3483   &11    &0.1355\\
100    &6    &0.1865    &2    &0.2856    &9    &0.1937\\
\hline 
\end{tabular}
\end{table}

\subsection{A linear algebraic system related to an inverse problem for the heat equation}

In this section, we apply the DSM and the VR to solve a linear algebraic
system used in the test problem {\it heat} from {\it Regularization
tools} in \cite{Hansen1}.  This linear algebraic system is a part of
numerical solutions to an inverse problem for the heat equation.  This
problem is reduced to a Volterra integral equation of the first kind with
$[0,1]$ as the integration interval.  The kernel is $K(s,t)=k(s-t)$ with
$$ k(t)=\frac{t^{-3/2}}{2\kappa \sqrt{\pi}}\exp(-\frac{1}{4\kappa^2 t}).
$$ Here, we use the default value $\kappa=1$. In this test in
\cite{Hansen1} the integral equation is discretized by means of simple
collocation and the midpoint rule with $n$ points.  
The unique exact solution
$u_n$ is constructed, and then the right-hand side $b_n$ is produced as
$b_n=A_nu_n$ (see \cite{Hansen1}).  In our test, we use $n=10,20,...,100$
and $b_{n,\delta} = b_n + e_n$, where $e_n$ is a vector containing random
entries, normally distributed with mean 0, variance 1, and scaled so that
$\|e_n\|=\delta_{rel}\|b_n\|$.  This linear system is ill-posed:  
the condition number of $A_{100}$ obtained by using the function
{\it cond} provided in MATLAB is $1.3717\times 10^{37}$.  As we have
discussed earlier, this condition number may be not accurate because
of the limitations of the program {\it cond} provided in MATLAB.  However, 
this number shows that the corresponding linear
algebraic system is ill-conditioned.

\begin{table}[h] 
\caption{Numerical results for inverse heat equation with $\delta_{rel}=0.05$, $n=10i, i=\overline{1,10}$.}
\label{heattab1}
\centering
\small
\begin{tabular}{@{  }c@{\hspace{2mm}}
@{\hspace{2mm}}|c@{\hspace{2mm}}c@{\hspace{2mm}}|c@{\hspace{2mm}}c@{\hspace{2mm}}|
c@{\hspace{2mm}}c@{\hspace{2mm}}|c@{\hspace{2mm}}r@{\hspace{2mm}}l@{}} 
\hline
&\multicolumn{2}{c|}{DSM}&\multicolumn{2}{c|}{VR$_{i}$}&\multicolumn{2}{c|}{VR$_{n}$}\\
$n$&
$N_{\text{linsol}}$&$\frac{\|u_\delta-y\|_{2}}{\|y\|_2}$&
$N_{\text{linsol}}$&$\frac{\|u_\delta-y\|_{2}}{\|y\|_2}$&
$N_{\text{linsol}}$&$\frac{\|u_\delta-y\|_{2}}{\|y\|_2}$\\
\hline
10    &8    &0.2051    &1    &0.2566    &6    &0.2066\\
20    &4    &0.2198    &1    &0.4293    &8    &0.2228\\
30    &7    &0.3691    &1    &0.4921    &6    &0.3734\\
40    &4    &0.2946    &1    &0.4694    &8    &0.2983\\
50    &4    &0.2869    &1    &0.4780    &7    &0.3011\\
60    &4    &0.2702    &1    &0.4903    &9    &0.2807\\
70    &4    &0.2955    &1    &0.4981    &6    &0.3020\\
80    &5    &0.2605    &1    &0.4743   &10    &0.2513\\
90    &5    &0.2616    &1    &0.4802    &8    &0.2692\\
100    &5    &0.2588    &1    &0.4959    &6    &0.2757\\
\hline 
\end{tabular}
\end{table}

Looking at the Table~\ref{heattab1} one can see that in some situations 
the VR$_{n}$ is not as accurate as the DSM even when it takes more iterations than the DSM. 
Overall, the results obtained by the DSM are often slightly more accurate than 
those by the VR$_{n}$. 
The time of computation of the DSM is comparable to that of the VR$_{n}$. 
In some situations, the results by VR$_{n}$ and the VR$_{i}$ are the same although it uses 3 more iterations than does the DSM.
The conclusion from this Table is that DSM competes favorably with the VR$_{n}$ in both accuracy and time of computation.

\begin{figure}[!h!tb]
\centerline{%
\includegraphics[scale=0.45]{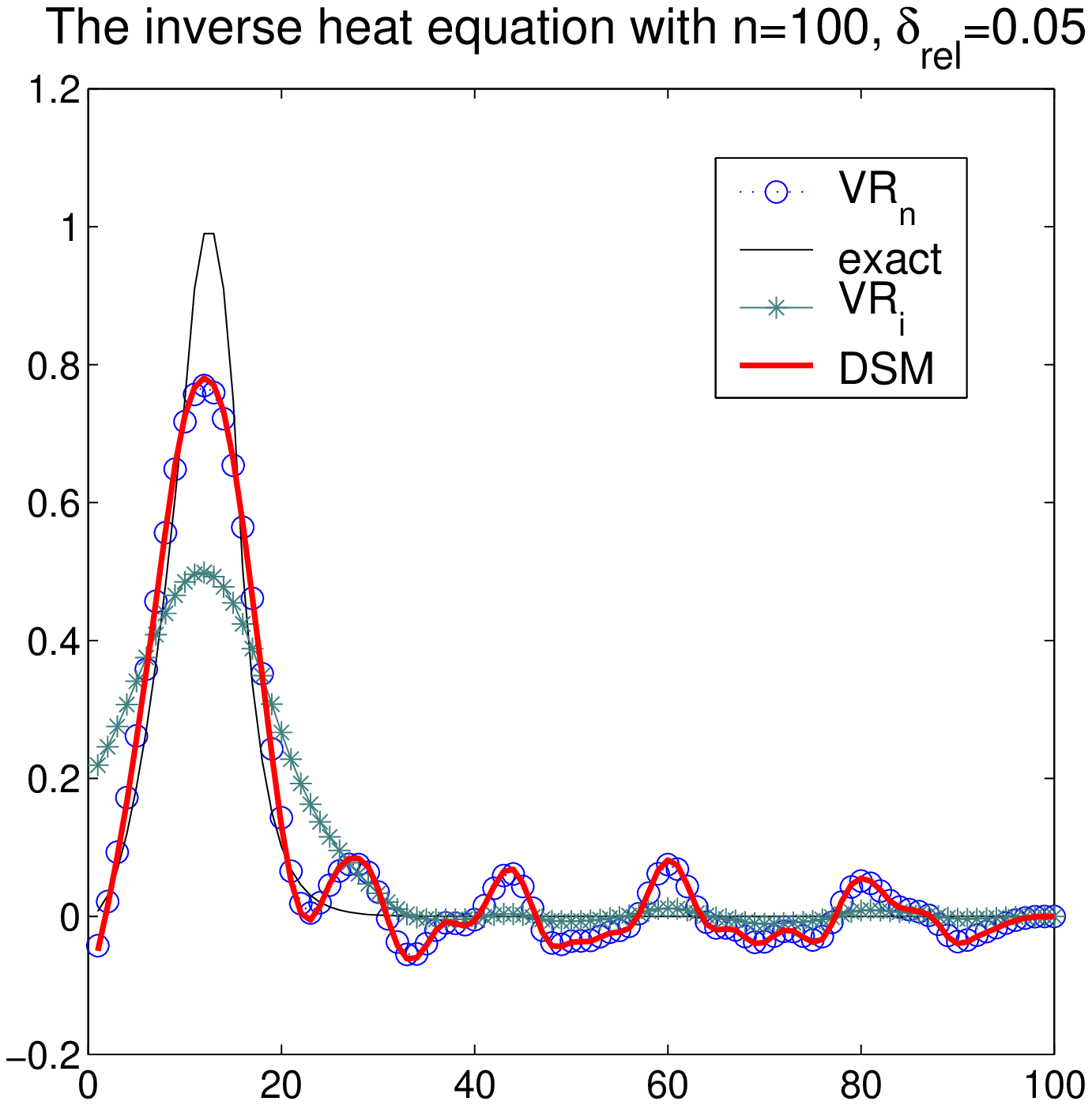}\quad\includegraphics[scale=0.45]{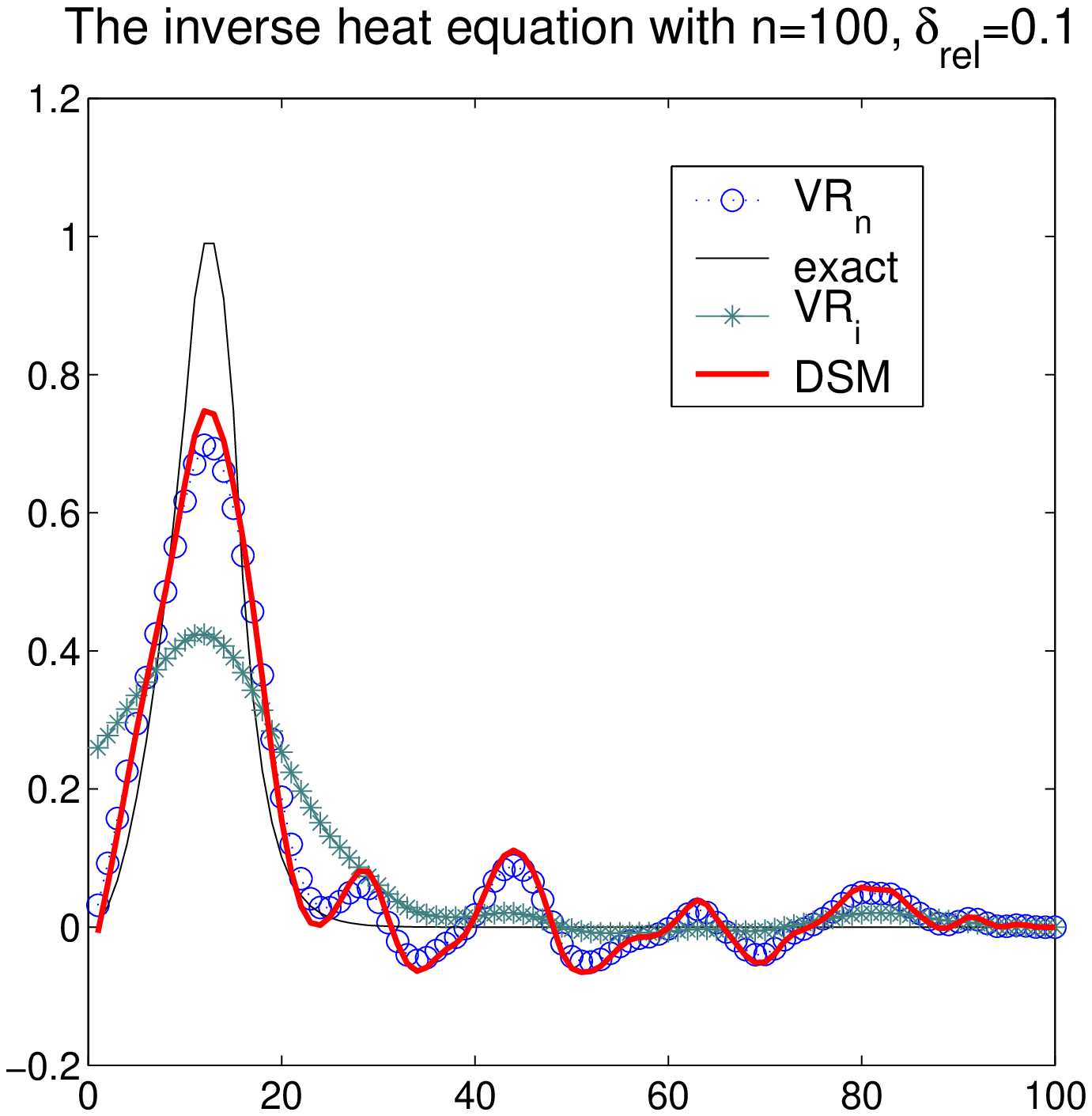}}
\caption{Plots of solutions obtained by DSM, VR for the inverse heat equation when $n=100$, $\delta_{rel}=0.05$ (left) and $\delta_{rel}=0.1$ (right).}
\label{figheat}
\end{figure}

Figure~\ref{figheat} plots numerical solutions to the inverse heat equation for $\delta_{rel}=0.05$ and $\delta_{rel}=0.02$ when $n=100$. From the figure we can see that the numerical solutions obtained by the DSM are about the same those by the VR$_{n}$. In these examples, the time of computation of the DSM is about the same as that of the VR$_{n}$.

The conclusion is that the DSM competes favorably with the VR$_{n}$ in this experiment.

\subsection{A linear algebraic system for the computation of the second derivatives}

Let us do some numerical experiments with linear algebraic systems arising in a
numerical experiment of computing the second derivative of a noisy function.

The problem is reduced to an integral equation of the first kind.
A linear algebraic system is obtained by a discretization of the integral equation whose kernel $K$ is Green's function
$$
K(s,t)=
\left\{
\begin{matrix}
s(t-1),\quad \text{if}\quad s<t\\
t(s-1),\quad \text{if}\quad s\geq t
\end{matrix}
\right. .
$$
Here $s,t\in[0,1]$ and as the right-hand side $f$ and the corresponding 
solution $u$ one chooses one of the following 
(see \cite{Hansen1}):
\begin{align*}
\text{case 1},\quad &f(s)=(s^3-s)/6,\quad u(t)=t,\\
\text{case 2},\quad &f(s)=e^s+(1-e)s-1,\quad u(t)=e^t,\\
\text{case 3},\quad &f(s)=
\left\{
\begin{matrix}
(4s^3-3s)/24,         &\quad \text{if} \quad&s<\frac{1}{2}\\
(-4s^3+12s^2-9s+1)/24,&\quad \text{if} \quad&s\geq\frac{1}{2}
\end{matrix}
\right .,\\
&u(t)=
\left\{
\begin{matrix}
t,&\quad \text{if} \quad&t<\frac{1}{2}\\
1-t,&\quad \text{if} \quad& t\geq\frac{1}{2}\\
\end{matrix}
\right..
\end{align*}

%Actually, in the linear equation $A_nu=f$, used in \cite{Hansen1}, the right-hand side $f$ contains a noise. 
%Although the noise level is small, it can happen that the error is underestimated or overestimated.
%Thus, we only use the matrix $A_n$ and $u$ and then we compute $\tilde{f}=A_nu$ as our linear equation.
Using $A_n$ and $u_n$ in \cite{Hansen1}, the right-hand side $b_n=A_nu_n$ 
is computed. 
Again, we use $n=10,20,...,100$ and $b_{n,\delta} = b_n + e_n$, 
where $e_n$ is a vector containing random entries, normally distributed with mean 0, variance 1, 
and scaled so that $\|e_n\|=\delta_{rel}\|b_n\|$. 
This linear algebraic system is mildly ill-posed: the condition 
number of $A_{100}$ is $1.2158\times10^{4}$.

Numerical results for the third case is presented in Table~\ref{case3}. 
In this case, the results obtained by the VR$_{n}$ are often slightly more accurate than those of the DSM. 
However, the difference between accuracy as well as the difference between time
of computation of these methods is small.
%The conclusion in this case is that the DSM is comparable with the VR$_{n}$.
Numerical results obtained by these two methods are much better than those of the VR$_{i}$. 

\begin{table}[h] 
\caption{Results for the deriv2 problem with $\delta_{rel}=0.01$, $n=100$ case 3.}
\label{case3}
\centering
\small
\begin{tabular}{@{  }c@{\hspace{2mm}}|c@{\hspace{2mm}}
c@{\hspace{2mm}}|c@{\hspace{2mm}}c@{\hspace{2mm}}
|c@{\hspace{2mm}}c@{\hspace{2mm}}|c@{\hspace{2mm}}r@{\hspace{2mm}}l@{}} 
\hline
&\multicolumn{2}{c|}{DSM}&\multicolumn{2}{c|}{VR$_{i}$}&\multicolumn{2}{c|}{VR$_{n}$}\\
$n$&
$N_{\text{linsol}}$&$\frac{\|u_\delta-y\|_{2}}{\|y\|_2}$&
$N_{\text{linsol}}$&$\frac{\|u_\delta-y\|_{2}}{\|y\|_2}$&
$N_{\text{linsol}}$&$\frac{\|u_\delta-y\|_{2}}{\|y\|_2}$\\
\hline
10    &4    &0.0500    &2    &0.0542    &6    &0.0444\\
20    &4    &0.0584    &2    &0.0708    &6    &0.0561\\
30    &4    &0.0690    &2    &0.0718    &6    &0.0661\\
40    &4    &0.0367    &1    &0.0454    &4    &0.0384\\
50    &3    &0.0564    &1    &0.0565    &4    &0.0562\\
60    &4    &0.0426    &1    &0.0452    &4    &0.0407\\
70    &5    &0.0499    &1    &0.0422    &5    &0.0372\\
80    &4    &0.0523    &1    &0.0516    &4    &0.0498\\
90    &4    &0.0446    &1    &0.0493    &4    &0.0456\\
100   &4    &0.0399    &1    &0.0415    &5    &0.0391\\
\hline 
\end{tabular}
\end{table}

\begin{figure}[!h!tb]
\centerline{%
\includegraphics[scale=0.45]{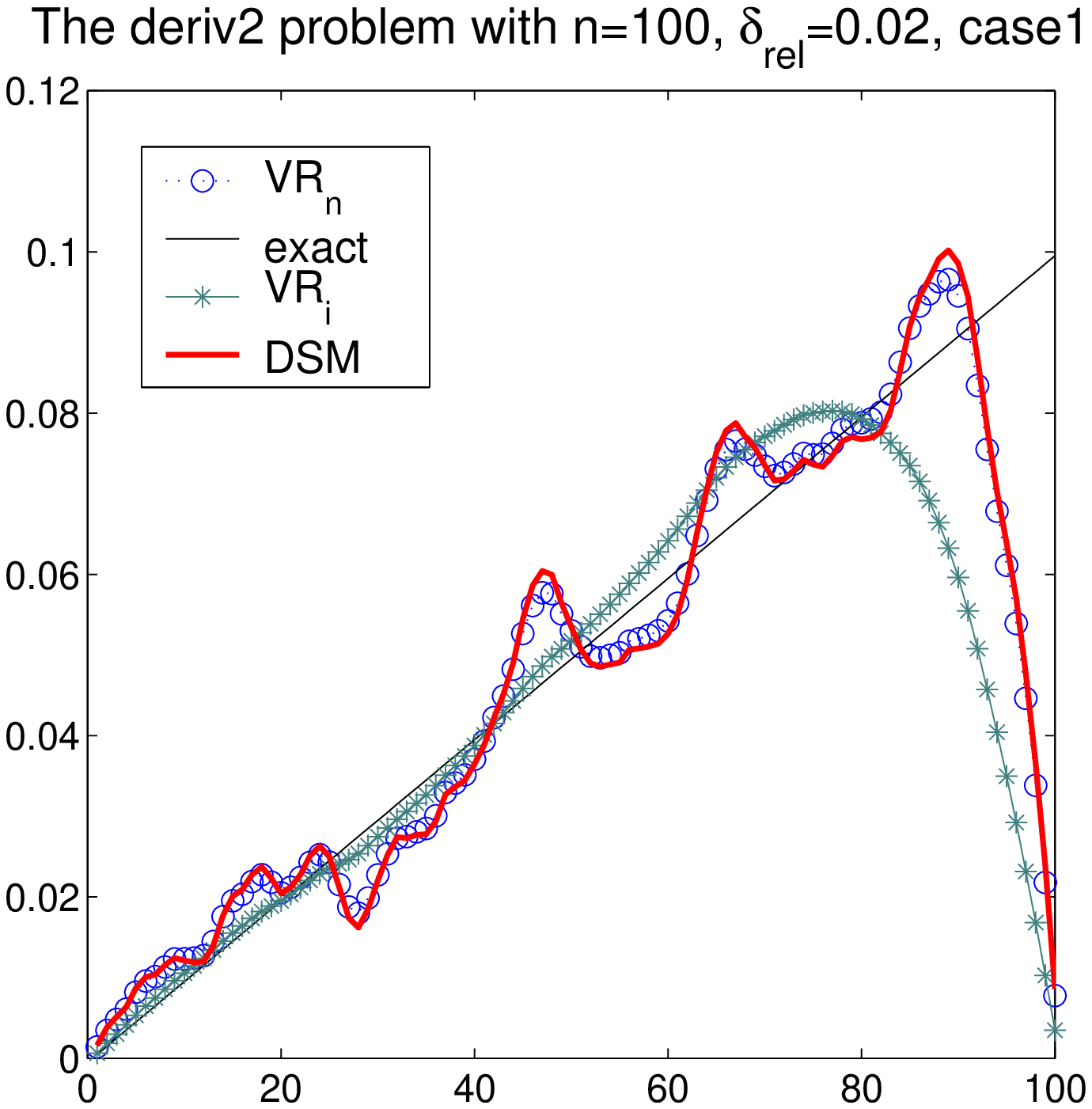}\quad\includegraphics[scale=0.45]{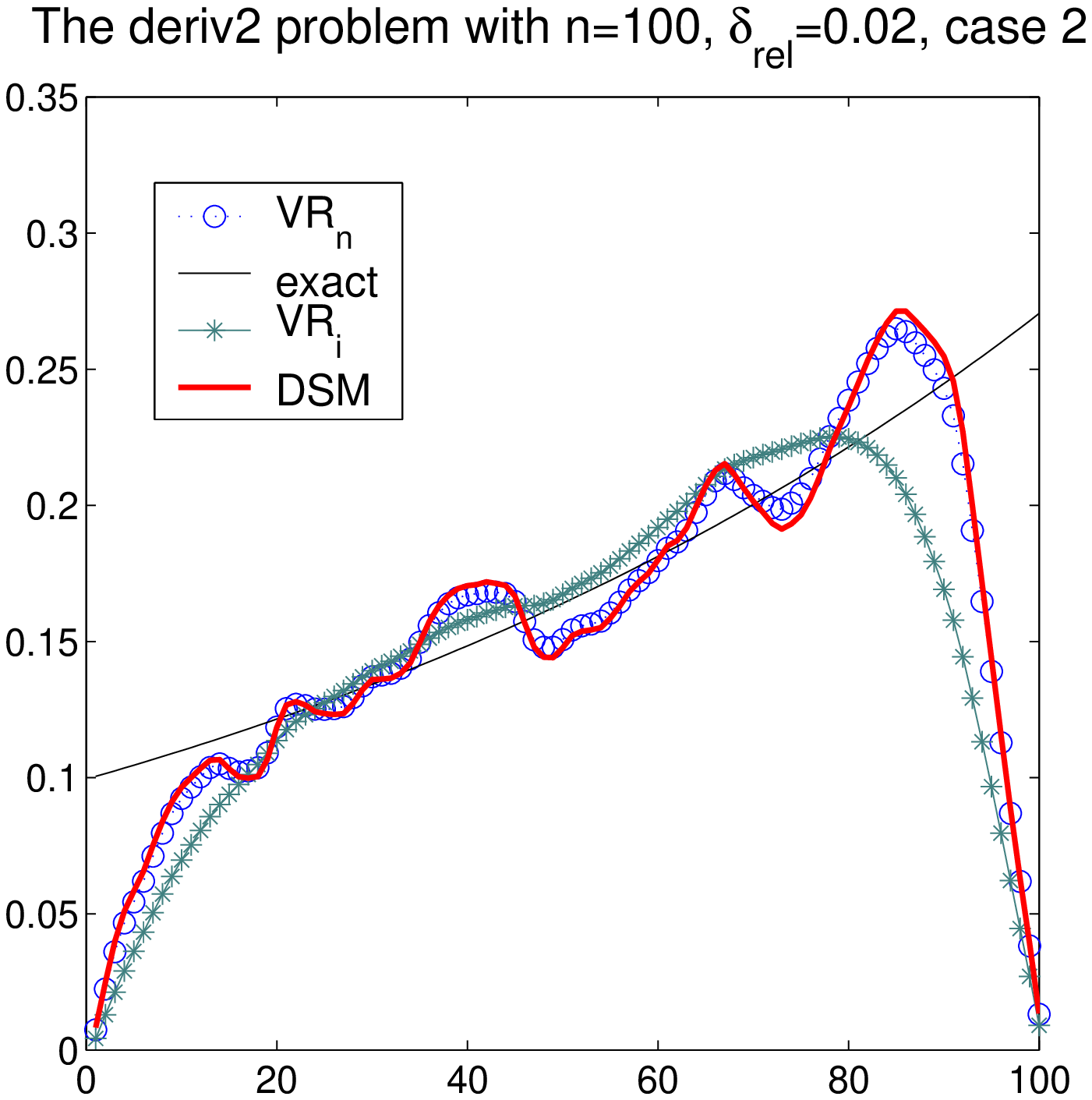}}
\caption{Plots of solutions obtained by DSM, VR for the deriv2 problem when $n=100$, 
$\delta_{rel}=0.02$ (left) and $\delta_{rel}=0.02$ (right).}
\label{figderiv}
\end{figure}

For other cases, case 1 and case 2, numerical results obtained by the DSM are slightly more accurate than those by the VR$_{i}$. 
Figure~\ref{figderiv}  plots the numerical solutions for these cases. 
%In case 2, the results by the DSM are more accurate than by the VR$_{n}$. 
The computation time of the DSM in these cases is about the same as or less than that of the VR$_{n}$.

The conclusion in this experiment is that the DSM competes favorably with the VR.
Indeed, the VR$_n$ is slightly better than the DSM in case 3 but slightly worse than the
DSM in cases 1 and 2.

\section{Concluding remarks}
The conclusions from the above experiments are:
\begin{enumerate}
\item{
The DSM always converges for $a(t)=\frac{a_0}{1+t}$ given that $a_0>a_{op}$. 
However, if $a_0$ is not well chosen, then the convergence speed may be  slow.
The parameter $a_0$ should be chosen so that it is greater than and close to 
the 
optimal $a_{op}$, i.e., the value minimizing the quantity:
\begin{equation*}
\|u_{a_{op}}-y\| = \inf_{a}\|u_a -y\|,\quad
u_a = T_a^{-1}Af_\delta.
\end{equation*}
However, since $a_{op}$ is not known and $a_M$ is often 
close to $a_{op}$, we choose $a_0$ so that
$$
\delta<\|AT_{a_0}^{-1}A^*f_\delta -f_\delta\|<2\delta.
$$
}
\item{
The DSM is sometimes faster than the VR. 
In general, the DSM is comparable with the VR$_{n}$ with respect to
computation time.} 
\item{
The DSM is often slightly more accurate than the VR, especially 
when $\delta$ is large. 
Starting with $a_0$ such that $\delta <\|AT_{a_0}^{-1}A^*f_\delta -f_\delta\|<2\delta$, 
the DSM often requires 4 to 7 iterations, and main cost in each iteration 
consists of solving the linear system $T_a u = A^*f_\delta$. 
The cost of these iterations is often about the same as the cost of using 
Newton's method to solve $a_M$ in the VR$_{n}$.}
\item{
For any initial $a_0$ such that 
$\delta <\|AT_{a_0}^{-1}A^*f_\delta -f_\delta\|<2\delta$, 
the DSM always converges to a solution which is often more accurate 
than that of the VR$_{n}$. 
However, with the same initial $a_0$, the VR$_{n}$ does not necessarily converge. 
In this case, we restart the Newton scheme to solve for the regularization parameter
with initial guess $a_1=\frac{a_0}{2}$ instead of $a_0$.}
%\item{
%Finding $a_0$ satisfying $\delta <\|A(A^*A+a_0 I)^{-1}A^*f_\delta -f_\delta\|
%<2\delta$ is easier 
%than finding $a_0$ such that the Newton scheme converges, 
%especially when $\delta_{rel}$ is large.}

\end{enumerate}

\label{lastpage}

\end{document}